\documentclass[reqno]{amsart}
\usepackage{setspace}
\usepackage{graphicx}
\usepackage{amsmath}
\usepackage{amsfonts}
\usepackage{amssymb}
\usepackage{amsthm}
\usepackage{mathtools}
\usepackage{subfigure}
\usepackage{enumerate}
\usepackage{enumitem}

\begin{document}

\centerline{\Huge\bf The Three and Fourfold Translative Tiles}

\medskip
\centerline{\Huge\bf  in Three-Dimensional Space}

\bigskip\medskip
\centerline{\Large Mei Han, Kirati Sriamorn, Qi Yang and Chuanming Zong}

\bigskip\bigskip
\centerline{\begin{minipage}{13cm}
{\bf Abstract.} This paper proves the following statement: {\it If a convex body can form a three or fourfold translative tiling in $\mathbb{E}^3$, it must be a parallelohedron.} In other words, it must be a parallelotope, a hexagonal prism, a rhombic dodecahedron, an elongated dodecahedron, or a truncated octahedron.
\end{minipage}}



\vspace{0.6cm}
\noindent
{\LARGE\bf 1. Introduction}

\bigskip\medskip
\noindent
Let $K$ be a convex body with interior $int(K)$ and boundary $\partial (K)$, and let $X$ be a discrete set, both in $\mathbb{E}^n$. We call $K+X$ a {\it translative tiling} of $\mathbb{E}^n$ and call $K$ a {\it translative tile} if $K+X=\mathbb{E}^n$ and the translates ${\rm int}(K)+{\bf x}_i$ are pairwise disjoint. In other words, if $K+X$ is both a packing in $\mathbb{E}^n$ and a covering of $\mathbb{E}^n$. In particular, we call $K+\Lambda$ a {\it lattice tiling} of $\mathbb{E}^n$ and call $K$ a {\it lattice tile} if $\Lambda $ is an $n$-dimensional lattice. Apparently, a translative tile must be a convex polytope. Usually, a lattice tile is called a {\it parallelohedron}.

In 1885, Fedorov \cite{fedo} characterized the two- and three-dimensional lattice tiles: {\it A two-dimensional lattice tile is either a parallelogram or a centrally symmetric hexagon; A three-dimensional lattice tile must be a parallelotope, a hexagonal prism, a rhombic dodecahedron, an elongated dodecahedron, or a truncated octahedron.}

The situations in higher dimensions turn out to be very complicated. Through the works of Delone \cite{delo}, $\check{S}$togrin \cite{stog} and Engel \cite{enge}, we know that there are exact $52$ combinatorially different types of parallelohedra in $\mathbb{E}^4$. A computer classification for the five-dimensional parallelohedra was announced by Dutour Sikiri$\acute{\rm c}$, Garber, Sch$\ddot{\rm u}$rmann and Waldmann \cite{dgsw} only in 2015.

Let $\Lambda $ be an $n$-dimensional lattice. The {\it Dirichlet-Voronoi cell} of $\Lambda $ is defined by
$$D=\left\{ {\bf x}: {\bf x}\in \mathbb{E}^n,\ \| {\bf x}, {\bf o}\|\le \| {\bf x}, \Lambda \|\right\},$$
where $\| X, Y\|$ denotes the Euclidean distance between $X$ and $Y$. Clearly, $D+\Lambda $ is a lattice tiling and the Dirichlet-Voronoi cell $D$ is a parallelohedron. In 1908, Voronoi \cite{voro} made a conjecture that {\it every parallelohedron is a linear image of the Dirichlet-Voronoi cell of a suitable lattice.} In $\mathbb{E}^2$, $\mathbb{E}^3$ and $\mathbb{E}^4$, this conjecture was confirmed by Delone \cite{delo} in 1929. In higher dimensions, it is still open.

To characterize the translative tiles is another fascinating problem. First it was shown by Minkowski \cite{mink} in 1897 that {\it every translative tile must be centrally symmetric}. In 1954, Venkov \cite{venk} and Aleksandrov \cite{alek} proved that {\it all translative tiles are parallelohedra.} Later, a new proof for this beautiful result was independently discovered by McMullen \cite{mcmu} (see also \cite{zong96}).

Let $X$ be a discrete multiset in $\mathbb{E}^n$ and let $k$ be a positive integer. We call $K+X$ a {\it $k$-fold translative tiling} of $\mathbb{E}^n$ and call $K$ a {\it $k$-fold translative tile} if every point ${\bf x}\in \mathbb{E}^n$ belongs to at least $k$ translates of $K$ in $K+X$ and every point ${\bf x}\in \mathbb{E}^n$ belongs to at most $k$ translates of ${\rm int}(K)$ in ${\rm int}(K)+X$. In other words, if $K+X$ is both a $k$-fold packing in $\mathbb{E}^n$ and a $k$-fold covering of $\mathbb{E}^n$. In particular, we call $K+\Lambda$ a {$k$-fold lattice tiling} of $\mathbb{E}^n$ and call $K$ a {\it $k$-fold lattice tile} if $\Lambda $ is an $n$-dimensional lattice. Apparently, a $k$-fold translative tile must be a convex polytope.

Multiple tilings were first investigated by Furtw\"angler \cite{furt} in 1936 as a generalization of Minkowski's conjecture on cube tilings. Let $C$ denote the $n$-dimensional unit cube. Furtw\"angler made a conjecture that {\it every $k$-fold lattice tiling $C+\Lambda$ has twin cubes. In other words, every multiple lattice tiling $C+\Lambda$ has two cubes which share a whole facet.} In the same paper, he proved the two- and three-dimensional cases. Unfortunately, when $n\ge 4$, this beautiful conjecture was disproved by Haj\'os \cite{hajo} in 1941. In 1979, Robinson \cite{robi} determined all the integer pairs $\{ n,k\}$ for which Furtw\"angler's conjecture is false. We refer to Zong \cite{zong05,zong06} for detailed accounts on this fascinating problem.

\medskip
Clearly, one of the most important and natural problems in multiple tilings is, for given $n$ and $k$, to classify or characterize all the $n$-dimensional $k$-fold lattice tiles and all the $n$-dimensional $k$-fold translative tiles (see Problems 1-4 at the end of Gravin, Robins and Shiryaev \cite{grs}). In the plane, it was proved by Yang and Zong \cite{yz1,yz2,zong20,zong-x} that, {\it besides parallelograms and centrally symmetric hexagons, there is no other two-, three- or fourfold translative tile. However, there are three classes of other fivefold translative tiles and three classes of other sixfold lattice tiles.} For more related results on multiple tilings, we refer to \cite{boll,gkrs,Grepstad,kolo,Lev-Liu}.

As a counterpart of the above result and as a generalization of the theorem of Venkov, Aleksandrov and McMullen, it was shown by Han, Sriamorn, Yang and Zong \cite{hsyz} that {\it a convex body can form a twofold translative tiling in $\mathbb{E}^3$ if and only if it is a parallelohedron}. This paper proves the following results.

\medskip\noindent
{\bf Theorem 1.} {\it In $\mathbb{E}^3$, every threefold translative tile is a parallelohedron. In other words, every three-dimensional threefold translative tile must be a parallelotope, a hexagonal prism, a rhombic dodecahedron, an elongated dodecahedron, or a truncated octahedron.}

\medskip\noindent
{\bf Theorem 2.} {\it In $\mathbb{E}^3$, every fourfold translative tile is a parallelohedron. In other words, every three-dimensional fourfold translative tile must be a parallelotope, a hexagonal prism, a rhombic dodecahedron, an elongated dodecahedron, or a truncated octahedron.}

\medskip\noindent
{\bf Theorem 3.} {\it Every belt of a fivefold translative tile in $\mathbb{E}^3$ has at most ten facets. In other words, every belt of a fivefold translative tile in $\mathbb{E}^3$ has four, six, eight or ten facets.}

\medskip\noindent
{\bf Theorem 4.} {\it Every belt of a sixfold translative tile in $\mathbb{E}^3$ has at most fourteen facets. In other words, every belt of a sixfold translative tile in $\mathbb{E}^3$ has four, six, eight, ten, twelve or fourteen facets.}

\vspace{0.6cm}
\noindent
{\Large\bf 2. Some Basic Known Results}

\medskip
\noindent
In this section, we recall some known concepts and results which will be useful for this paper.

\medskip
In 1885, E. S. Fedorov studied the two- and three-dimensional lattice tiles. He proved the following result.

\smallskip\noindent
{\bf Lemma 1 (Fedorov \cite{fedo}).} {\it A two-dimensional lattice tile is either a parallelogram or a centrally symmetric hexagon; A three-dimensional lattice tile must be a parallelotope, a hexagonal prism, a rhombic dodecahedron, an elongated dodecahedron, or a truncated octahedron.}

\smallskip
Tilings in higher dimensions have been studied by Minkowski \cite{mink}, Voronoi \cite{voro}, Delone \cite{delo}, Venkov \cite{venk}, Alexsandrov \cite{alek}, McMullen \cite{mcmu} and many others. This paper needs the following concepts and results of them.

\smallskip
\noindent
{\bf Definition 1.} Let $P$ denote an $n$-dimensional centrally symmetric convex polytope with centrally symmetric facets and let $V$ denote a $(n-2)$-dimensional face of $P$. We call the collection of all those facets of $P$ which contain a translate of $V$ as a subface a belt of $P$.

\medskip\noindent
{\bf Lemma 2 (Venkov \cite{venk} and McMullen \cite{mcmu}).} {\it  Let $K$ be an $n$-dimensional convex body. The following three statements are equivalent to each other:}
\begin{enumerate}
\item {\it $K$ is a translative tile;}
\item {\it $K$ is a centrally symmetric polytope with centrally symmetric facets such that each belt contains four or six facets;}
\item {\it $K$ is a parallelohedron.}
\end{enumerate}

\medskip\noindent
{\bf Definition 2.} Let $P$ be an $n$-dimensional convex polytope. We call it a zonotope if it is a Minkowski sum of finite number of segments.
In other words,
$$P=S_1+S_2+\ldots +S_w,$$
where $w$ is an integer and $S_i$ are segments in $\mathbb{E}^n$.

\medskip
In 2012, N. Gravin, S. Robins and D. Shiryaev studied multiple tilings in general dimensions and discovered the following result.

\medskip\noindent
{\bf Lemma 3 (Gravin, Robins and Shiryaev \cite{grs}).}  {\it An $n$-dimensional $k$-fold translative tile is a centrally symmetric polytope with centrally symmetric facets. In particular, in $\mathbb{E}^3$, every $k$-fold translative tile is a zonotope.}

\vspace{0.6cm}
\noindent
{\Large\bf 3. Some Properties of Zonotopes}

\medskip
\noindent
In this section, we introduce some basic properties of zonotopes which will be useful for our proofs. Let $P$ be a zonotope in $\mathbb{E}^3$. By a translation, without loss of generality, we may assume that
$$P=\left\{\sum_{i=1}^wx_i{\bf u}_i, \quad x_i\in [0,1]\right\},\eqno(1)$$
where $U=\{{\bf u}_1, {\bf u}_2, \ldots , {\bf u}_w\}$ is a set of pairwise linear independent vectors in $\mathbb{E}^3$ and each edge of $P$ is a translate of exactly one generator ${\bf u}_i$. For convenience, we call $U$ a generator set of $P$.

\medskip\noindent
{\bf Lemma 4.} {\it Let ${\bf q}$ be a vertex of $P$. If
$${\bf q}=\sum_i x_i{\bf u}_i,\quad x_i\in[0,1],$$
then all $x_i$ are either $0$ or $1$.}

\medskip
\noindent
{\bf Proof.} If, on the contrary, there is an index $j$ such that $x_j\in(0,1)$. We take
$${\bf q}'=\sum_{i\neq j} x_i{\bf u}_i,$$
and
$${\bf q}^*={\bf u}_j+\sum_{i\neq j} x_i{\bf u}_i.$$
It is clear that both ${\bf q}'$ and ${\bf q}^*$ belong to $P$ and
$${\bf q}= (1-x_j){\bf q}'+x_j {\bf q}^*,\eqno(2)$$
which contradicts the assumption that ${\bf q}$ is a vertex of $P$. The lemma is proved.\hfill{$\Box$}

\bigskip
Let $G$ be an edge of $P$ and let $B(G)$ be the belt determined by $G$, which consists of $2m$ facets $F_1, F_2, \ldots, F_{2m}$ enumerated in clockwise order. Let $G_1, G_2, \ldots, G_{2m}$ be the translates of $G$ such that $G_i, G_{i+1} \subset F_i$ and $G_1 = G$. Let ${\bf q}_1$ and ${\bf q}'_1$ be the two endpoints of $G_1$. By Lemma 4, we may assume that
$${\bf q}_1={\bf u}_{i_1}+{\bf u}_{i_2}+\cdots+{\bf u}_{i_t}$$
and
$${\bf q}'_1={\bf u}_1+{\bf u}_{i_1}+{\bf u}_{i_2}+\cdots+{\bf u}_{i_t},$$
where ${\bf u}_1$ is the generator which parallel to $G$. Observe that $P-{\bf q}_1$ is a translate of $P$ with a generator set
$$(U\setminus\{{\bf u}_{i_1}, {\bf u}_{i_2}, \ldots, {\bf u}_{i_t}\})\cup \{-{\bf u}_{i_1}, -{\bf u}_{i_2}, \ldots, -{\bf u}_{i_t}\}.$$
Therefore, without loss of generality, we may assume that ${\bf q}_1$ is the origin and ${\bf q}'_1={\bf u}_1$.

Let $G$ be an edge of $P$ with an endpoint at origin, and let $B(G)$ be the belt determined by $G$, which consists of $2m$ facets $F_1, F_2, \ldots, F_{2m}$. Let $U=\{{\bf u}_1, {\bf u}_2, \ldots , {\bf u}_w\}$ be a generator set of $P$. For $1\leq i\leq m$, let ${\bf u}_0^i, {\bf u}_1^i, {\bf u}_2^i, \ldots, {\bf u}_{k_i}^i$ be the generators in $U$ which are parallel to the edges of $F_i$, where ${\bf u}_0^i$ is a translate of $G$. For convenience, we assume that ${\bf u}_0^i, {\bf u}_1^i, {\bf u}_2^i, \ldots, {\bf u}_{k_i}^i$ are in clockwise order and define
$$R_i=\left\{\sum_{j=1}^{k_i}x_j{\bf u}_j^i, \quad x_j\in[0,1]\right\}.\eqno(3)$$

\medskip\noindent
{\bf Lemma 5.} {\it For $2\leq i\leq m$, we have}
$$F_1+R_i\subset P.$$

\medskip\noindent
{\bf Proof.} Noting that
$$F'_1=\left\{\sum_{j=0}^{k_1}x_j{\bf u}_j^1, \quad x_j\in[0,1]\right\}\eqno(4)$$
is a subset of $P$ containing $G_1$. Furthermore, since ${\bf u}_0^1, {\bf u}_1^1, \ldots, {\bf u}_{k_1}^1$ are the generators of $P$ which are parallel to the edges of $F_1$, we have that $F'_1=F_1$. When $2\leq i\leq m$, it is easy to see that  ${\bf u}_0^1,$ ${\bf u}_1^1,$ ${\bf u}_2^1,$ $\ldots,$ ${\bf u}_{k_1}^1,$ ${\bf u}_1^i,$ ${\bf u}_2^i,$ $\ldots,$ ${\bf u}_{k_i}^i$ are pairwise linear independent. Consequently, we have
$$\left\{ {\bf u}_0^1, {\bf u}_1^1, {\bf u}_2^1, \ldots, {\bf u}_{k_1}^1, {\bf u}_1^i, {\bf u}_2^i, \ldots, {\bf u}_{k_i}^i\right\}\subset \left\{ {\bf u}_1, {\bf u}_2, \ldots, {\bf u}_w\right\}$$
and
$$F_1+\left\{\sum_{j=1}^{k_i}x_j{\bf u}_j^i, \quad x_j\in[0,1]\right\}\subset \left\{\sum_{i=1}^wx_i{\bf u}_i,\quad x_i\in[0,1]\right\}=P,\eqno(5)$$
as shown in Figure 1. The lemma is proved. \hfill{$\Box$}

\begin{figure}[h!]
\centering
\includegraphics[scale=.45]{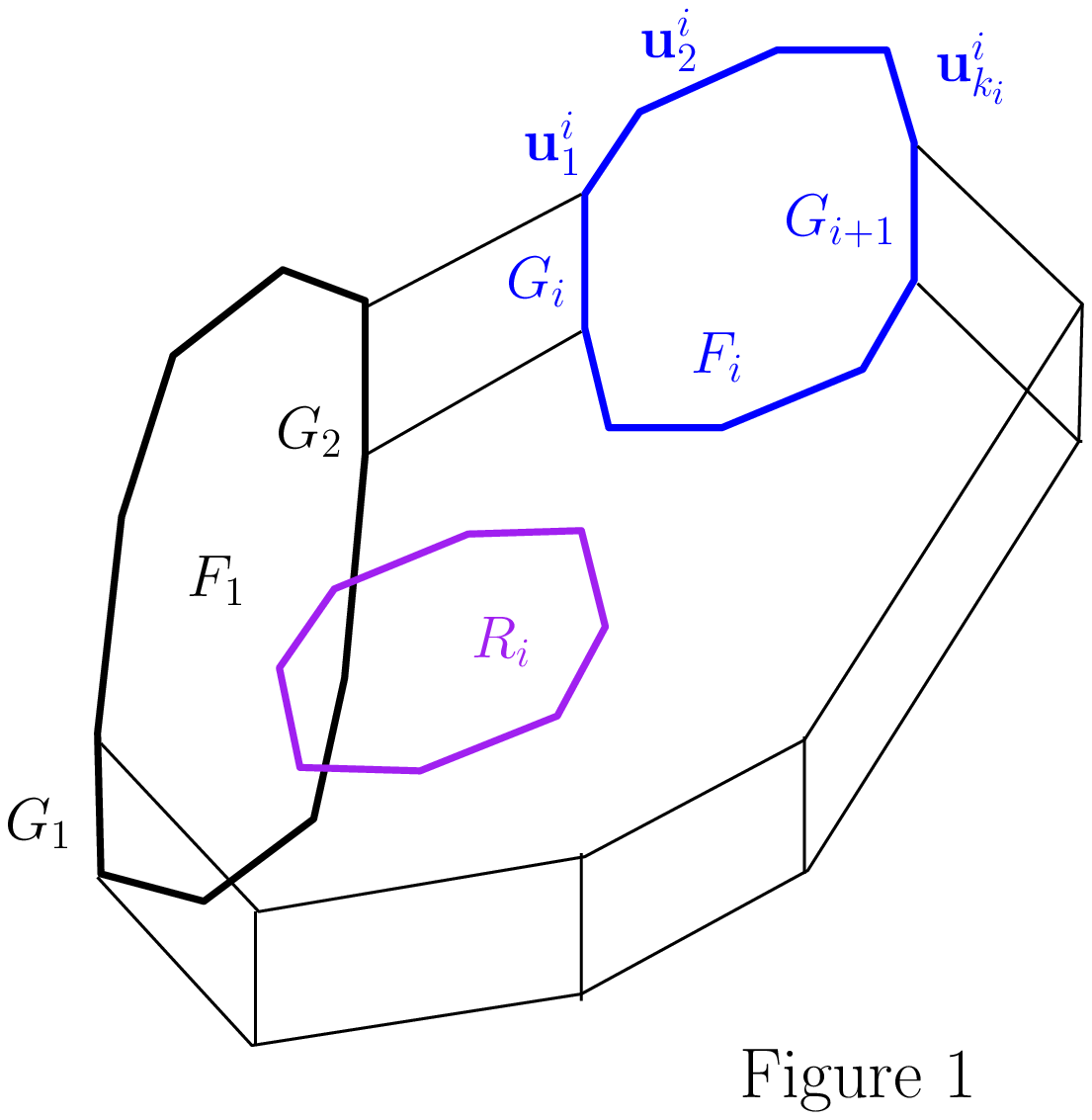}
\end{figure}

\medskip
For convenience, let $rint (F_i)$ denote the relative interior of $F_i$ and let $R'_i$ denote the set $R_i\setminus \{{\bf o}\}$.

\medskip\noindent
{\bf Corollary 1.} {\it If $m\ge 3$ and $2\le i\le m$, then}
$$rint (F_1)+R'_i\in int (P).$$
	
\medskip
For $1\leq i\leq m$, it is easy to see that ${\bf u}_0^i,$ ${\bf u}_1^i,$ $\ldots,$ ${\bf u}_{k_i}^i$ are also the generators in $U$ which are parallel to the edges of $F_{m+i}$. For convenience, we define ${\bf u}_{0}^{m+i}={\bf u}_0^i$ and ${\bf u}_{k}^{m+i}=-{\bf u}_{k_i+1-k}^i$, for $k=1,2,\ldots, k_i$. By symmetry, we take $R_{m+i}=-R_{i}$. Then, for $1\leq i\leq 2m$, one can deduced by induction that
$$F_i=\sum_{t=1}^{i-1}\sum_{s=1}^{k_t}{\bf u}_s^t +Q_i=G_i+R_i,\eqno(6)$$
where
$$Q_i=\left\{\sum_{j=0}^{k_i}x_j{\bf u}_j^i,\quad x_j\in[0,1]\right\}.\eqno(7)$$

For $1\leq i\leq 2m$, we abbreviate ${\bf u}_1^i+{\bf u}_2^i+\cdots+{\bf u}_{k_i}^i$ to ${\bf g}_i$. Then,
$$\sum_{j=1}^{i-1}{\bf g}_j+\sum_{j=0}^k{\bf u}_j^i, \quad 0\leq k\leq k_i $$
and
$$\sum_{j=1}^{i-1}{\bf g}_j+\sum_{j=k}^{k_i}{\bf u}_j^i, \quad 0\leq k\leq k_i$$
are the $2(k_i+1)$ vertices of $F_i$.

\medskip\noindent
{\bf Definition 3.} Let $F$ be a centrally symmetric polygon with an edge $E$ and let ${\bf g}$ be a vector such that $E+{\bf g}$ is also an edge of $F$. Then we call ${\bf g}$ the translation vector of $E$ in $F$. For ${\bf v}\in rint (E)$, we call ${\bf v}+{\bf g}$ the corresponding point of ${\bf v}$ in $F$.

\medskip
Let $F$ be a centrally symmetric convex $(2k+2)$-gon with edges $E_0,$ $E_1,$ $\ldots,$ $E_{2k+1}$ enumerated in clockwise order. Suppose that ${\bf v}_0$ is the common endpoint of edges $E_0$ and $E_{2k+1}$. Let ${\bf u}_0,$ ${\bf u}_1,$ $\ldots,$ ${\bf u}_k$ be the vectors which are parallel to $E_0,$ $E_1,$ $\ldots,$ $E_{k}$, respectively. One can easily deduce the following lemma which will be useful later. Its proof can be illustrated by Figure 2.

\begin{figure}[h!]
\centering
\includegraphics[scale=.45]{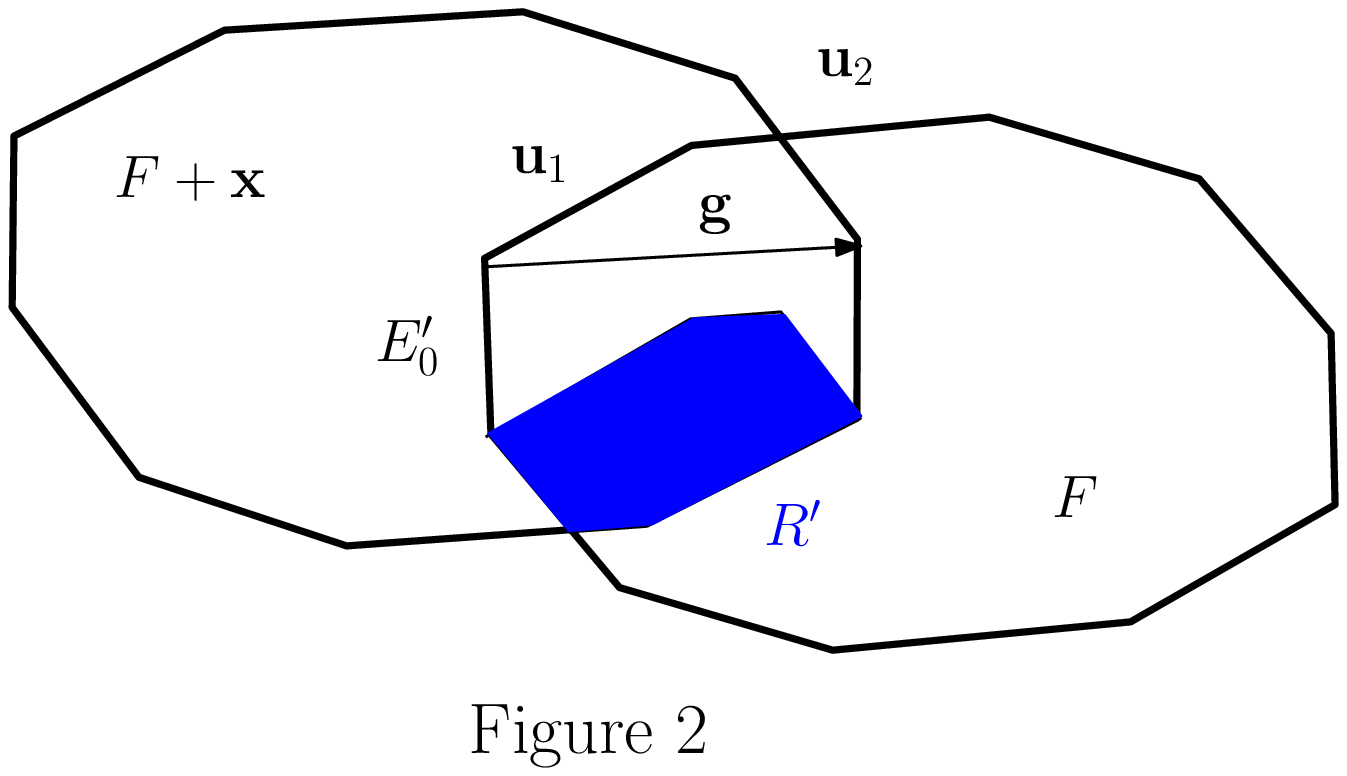}
\end{figure}

\medskip\noindent
{\bf Lemma 6.} {\it Suppose that ${\bf x}$ is a point such that $E_0\cap (F+{\bf x})\neq\emptyset$. We write $E'_0=E_0\cap (F+{\bf x})$. Then there exist $1\leq s<t\leq k$ and $0\leq\alpha<1,$ $0\leq\beta<1$ such that
$$F\cap(F+{\bf x})=E'_0+R',$$
where
$$R'=\left\{\left(\sum_{j=1}^{s-1}x_j{\bf u}_j\right)+x_s(\alpha {\bf u}_s)+x_t(\beta {\bf u}_t)+\left(\sum_{j=t+1}^{k}x_j{\bf u}_j\right), \quad  x_j\in[0,1]\right\}.$$
In particular, the vector
$${\bf g}=\left(\sum_{j=1}^{s-1}{\bf u}_j\right)+\alpha {\bf u}_s+\beta {\bf u}_t+\left(\sum_{j=t+1}^{k}{\bf u}_j\right)$$
is the translation vector of $E'_0$ in $F\cap(F+{\bf x})$ .}

\vspace{0.6cm}
\noindent
{\Large\bf 4. Adjacent Wheels and Proper Points}

\medskip\noindent
Assume that $G$ is parallel to the $z$-axis and ${\bf u}_0^1$ is in the positive $z$-direction. Let $K(G)$ denote the union of the facets which are not contained in $B(G)$ and let $X$ be a multiset such that $P+X$ is a $\tau$-fold tiling, for some positive integer $\tau$. For convenience, we define
$$\mathcal{G}=\{G_1, G_2, \ldots, G_{2m}\}$$
and
$$\mathcal{G}+X=\{G'+{\bf x}:\ G'\in\mathcal{G},\ {\bf x}\in X\}.$$
	
In the study of two-dimensional multiple tilings, the concept of adjacent wheel introduced by Yang and Zong \cite{yz2} played a crucial role. In $\mathbb{E}^3$, its analogy is also very useful. Assume that $G'\in \mathcal{G}+X$ and ${\bf v}\in G'$, where ${\bf v}\not\in K(G)+X$.  Let $X^{\bf v}$ denote the subset of $X$ consisting of all points ${\bf x}_i$ such that
$${\bf v}\in B(G)+{\bf x}_i.$$
Since $P+X$ is a multiple tiling, the set $X^{\bf v}$ can be divided into disjoint subsets $X_1^{\bf v},$ $X_2^{\bf v},$ $\ldots,$ $X_t^{\bf v}$ such that the translates in $P+X_j^{\bf v}$ can be re-enumerated as $P+{\bf x}_1^j,$ $P+{\bf x}_2^j,$ $\ldots,$ $P+{\bf x}_{s_j}^j$ satisfying the following conditions:
\begin{enumerate}
\item ${\bf v}\in B(G)+{\bf x}_i^j$ for all $i=1,2,\ldots,s_j$;
\item Let $\angle_i^j$ denote the inner dihedral angle of $P+{\bf x}_i^j$ at ${\bf v}$ with two half-planes $H_{i,1}^j$ and $H_{i,2}^j$ such that $H_{i,1}^j$ and $H_{i,2}^j$ are clockwise when viewed from above the positive $z$-axis. Then the angles join properly as
	$$H_{i,2}^j=H_{i+1,1}^j$$
	holds for all $i=1, 2, \ldots, s_j$, where $H_{s_j+1,1}^j=H_{1,1}^j$.
\end{enumerate}
We call such a sequence $P+{\bf x}_1^j,$ $P+{\bf x}_2^j,$ $\ldots,$ $P+{\bf x}_{s_j}^j$ an \textit{adjacent wheel} at ${\bf v}$. Then we define
$$\varpi({\bf v})=\frac{1}{2\pi}\sum_{j=1}^t\sum_{i=1}^{s_j}\angle_i^j\eqno(8)$$
and
$$\varphi({\bf v})=\sharp\{{\bf x}_i:\  {\bf x}_i\in X, {\bf v}\in int (P)+{\bf x}_i\}.\eqno(9)$$
Clearly, we have
$$\tau = \varpi({\bf v}) + \varphi({\bf v}).\eqno(10)$$

\medskip	
From Lemma 5 (Corollary 1), similar to Lemma 2.1 and Lemma 2.2 of Yang and Zong \cite{yz2}, one can obtain the following results.

\medskip
\noindent
{\bf Lemma 7.} {\it Assume that ${\bf v}\in G_i+X$, where $i\in\{1,\ldots,2m\}$ and ${\bf v}\not\in K(G)+X$. There are at least $\lceil (m-3)/2\rceil$ different translates $P+{\bf x}$ satisfying both
$${\bf v}\in B(G)+{\bf x}$$
and}
$$({\bf v}+R_i)\setminus\{{\bf v}\}\subset  int (P)+{\bf x}.$$

\medskip
\noindent
{\bf Lemma 8.} {\it Let $G'\in \mathcal{G}+X$ and ${\bf v}\in G'\setminus(K(G)+X)$. Then
$$\varpi({\bf v})=\kappa\cdot\frac{m-1}{2}+\ell\cdot\frac12,$$
where $\kappa$ is a positive integer and $\ell$ is a non-negative integer. In fact, $\ell $ is the number of the facets in $\mathcal{F}+X$ which take ${\bf v}$ as a relative interior point.}
	
\medskip
For convenience, we introduce the following definition.

\medskip\noindent
{\bf Definition 4.}	Suppose that ${\bf x}\in X$, ${\bf v}\in B(G)+{\bf x}$ and ${\bf v}\not\in K(G)+{\bf x}$. If ${\bf v}\in G'$ for some $G'\in \mathcal{G}+{\bf x}$, then we call $P+{\bf x}$ an \textit{E-type piece} at ${\bf v}$. If ${\bf v}\not\in G'$ for all $G'\in\mathcal{G}+{\bf x}$, then we call $P+{\bf x}$ a \textit{F-type piece} at ${\bf v}$.

\medskip	
Obviously, in Lemma 8, $\ell $ is the number of F-type pieces at ${\bf v}$.

\medskip
\noindent
{\bf Lemma 9.} {\it Assume that $G'\in\mathcal{G}+X$.  Most of the points ${\bf v}\in G'$ (with only finite number of exceptions) satisfy the following conditions:}
\begin{enumerate}
\item  ${\bf v}\not\in  K(G)+X$;
\item  {\it If ${\bf v}\in rint (F_i)+{\bf x}$ holds for some integer $i$ and some point ${\bf x}\in X$, then there is ${\bf v}^*\in (G_i+{\bf x})\setminus (K(G)+X)$ and ${\bf v}^{**}={\bf v}^*+{\bf g}_i\in (G_{i+1}+{\bf x})\setminus (K(G)+X)$ such that
    $${\bf v}\in {\bf v}^*+R_i={\bf v}^{**}-R_i;$$}
\item  {\it If ${\bf v}\in G^*\cap F\subset F^*\cap F$, where $F\in B(G)+X$, $F^*\in B(G)+X$, and $G^*$ is an edge of $F^*$ which is parallel to $G'$, then the corresponding point of ${\bf v}$ in $F^*\cap F$ is not in $K(G)+X$.}
\end{enumerate}

\medskip\noindent
{\bf Proof.} Since $P+X$ is a $\tau$-fold tiling where $\tau$ is a positive integer, there are only a finite number of ${\bf x}\in X$ such that $G'\cap (P+{\bf x})\neq\emptyset$.

By the definition of $K(G)$, we know that $G'$ is not parallel to any facets in $K(G)$. Hence, one obtains that $G'\cap (K(G)+X)$ is finite.
	
Suppose that $L=G'\cap (rint (F_i)+{\bf x})\neq\emptyset$. Since $F_{i}+{\bf x}=(G_i+{\bf x})+R_i$, there must be a segment $L^*\subset G_i+{\bf x}$ and a vector ${\bf g}\in R_i$ such that $L=L^*+{\bf g}$. For any ${\bf v}\in L$, if we take ${\bf v}^*={\bf v}-{\bf g}\in L^*\subset G_i+{\bf x}$ and ${\bf v}^{**}={\bf v}-{\bf g}+{\bf g}_i\in G_{i+1}+{\bf x}$, then
$${\bf v}={\bf v}^*+{\bf g}\in {\bf v}^*+R_i={\bf v}^{**}-R_i.\eqno(11)$$
Noting that both $(G_i+{\bf x})\cap (K(G)+X)$ and $(G_{i+1}+{\bf x})\cap (K(G)+X)$ have finite number of points. Furthermore, there are only finite numbers of integers $i$ and ${\bf x}\in X$ such that $L\neq\emptyset$. It follows that there are only a finite number of points ${\bf v}\in G'$ which do not satisfy the second condition.
	
Suppose that $L= G'\cap G^*\cap F\neq\emptyset$. Denote by $L^{**}$ the set of all corresponding points of ${\bf v}\in L$ in $F^*\cap F$. It is clear that $L^{**}\cap(K(G)+X)$ is finite. Moreover, there are only finite numbers of $F\in B(G)+X$ and $F^*\in B(G)+X$ such that $L\neq\emptyset$. Hence, one can deduce that there are only a finite number of points ${\bf v}\in G'$ which do not satisfy the third condition.

In conclusion, there are only finite points in $G'$ which do not satisfy the desired conditions. The lemma is proved.\hfill{$\Box$}

\medskip\noindent
{\bf Definition 5.} Assume that  $G'\in\mathcal{G}+X$ and ${\bf v}\in G'$. If ${\bf v}$ satisfies all of the conditions in Lemma 9, then we call it a \textit{proper point} on $G'$.

\medskip
\noindent
{\bf Lemma 10.} {\it Assume that $G'\in\mathcal{G}+X$. If ${\bf v}$ is a proper point on $G'$, then}
$$\varphi({\bf v})\geq\left\lceil{\frac{m-3}{2}}\right\rceil.$$

\medskip
\noindent
{\bf Proof.} Without loss of generality, we may assume that ${\bf v}\in G_1+{\bf x}$. Let ${\bf v}^*$ be the corresponding point of ${\bf v}$ in $F_1+{\bf x}$, i.e., ${\bf v}^*={\bf v}+{\bf g}_1\in G_2+{\bf x}$. Since ${\bf v}$ is a proper point on $G'$, one obtains that ${\bf v}^*\not\in K(G)+X$. By applying Lemma 7 to the point ${\bf v}^*$, we have that there are at least $\lceil{(m-3)/2}\rceil$ different translates $P+{\bf x}_i$ satisfying both
$${\bf v}^*\in B(G)+{\bf x}_i$$
and
$$({\bf v}^*-R_{1})\setminus\{{\bf v}^*\}\subset  int (P)+{\bf x}_i.\eqno(12)$$
Since ${\bf v}={\bf v}^*-{\bf g}_1\in ({\bf v}^*-R_{1})\setminus\{{\bf v}^*\}$, we get
$${\bf v}\in int (P)+{\bf x}_i.$$
Thus, we have
$$\varphi({\bf v})\geq \left\lceil{\frac{m-3}{2}}\right\rceil.\eqno(13)$$
The lemma is proved. \hfill{$\Box$}

\vspace{0.6cm}
\noindent
{\Large\bf 5. Proofs of the Theorems}

\medskip\noindent
Recall that $P+X$ is a $\tau$-fold tiling, $G$ is an edge of $P$, and $B(G)$ is the belt of $2m$ facets determined by $G$. We may assume, without loss of generality, that ${\bf o}\in X$ and $G$ is perpendicular to the $xy$-plane. By adapting the method of Yang and Zong \cite{yz2}, one obtains the following results.

\medskip
\noindent
{\bf Lemma 11.} {\it If $m=4$, then $\tau\geq 5$.}

\medskip\noindent
{\bf Proof.} Assume that $G'\in \mathcal{G}+X$ and let ${\bf v}$ be a proper point on $G'$. First of all, we recall that
$$\tau=\varpi({\bf v})+\varphi({\bf v}).\eqno(14)$$
Moreover, from Lemma 8, we have that
$$\varpi({\bf v})=\kappa\cdot\frac32+\ell\cdot\frac12\geq 2.\eqno(15)$$
Now, we consider the following three cases.
		
\medskip\noindent
{\bf Case 1.} {\it $\varpi({\bf v})\geq 4$ for some proper point ${\bf v}$.} By Lemma 10, we get
$$\varphi({\bf v})\geq\left\lceil{\frac{m-3}{2}}\right\rceil =1$$
and hence
$$\tau=\varpi({\bf v})+\varphi({\bf v})\geq 4+1=5.\eqno(16)$$
		
\medskip\noindent
{\bf Case 2.} {\it $\varpi({\bf v})=3$ for some proper point ${\bf v}$.} By Lemma 8, we have
$$\varpi({\bf v})=\kappa\cdot\frac{3}{2}+\ell\cdot\frac12 =3.\eqno(17)$$

If $\ell>0$, then there is an F-type piece $P+{\bf x}$ at ${\bf v}$, where ${\bf x}\in X$. Suppose that ${\bf v}\in rint (F_i)+{\bf x}$. Using Lemma 9, one can choose ${\bf v}^*\in (G_i+{\bf x})\setminus(K(G)+X)$ and ${\bf v}^{**}={\bf v}^*+{\bf g}_i\in (G_{i+1}+{\bf x})\setminus(K(G)+X)$ such that
$${\bf v}\in {\bf v}^*+R_i={\bf v}^{**}-R_i.\eqno(18)$$
By applying Lemma 7 to ${\bf v}^*$ and ${\bf v}^{**}$, respectively, there must be ${\bf x}^*, {\bf x}^{**}\in X$ such that
$${\bf v}^*\in B(G)+{\bf x}^*,$$
$${\bf v}^{**}\in B(G)+{\bf x}^{**},$$
$$({\bf v}^*+R_{i})\setminus\{{\bf v}^*\}\subset  int (P)+{\bf x}^*$$
and
$$({\bf v}^{**}-R_{i})\setminus\{{\bf v}^{**}\}\subset  int (P)+{\bf x}^{**}.$$
Since ${\bf v}\in {\bf v}^*+R_i={\bf v}^{**}-R_i$, we obtain ${\bf v}\in int (P)+{\bf x}^*$ and ${\bf v}\in int (P)+{\bf x}^{**}$. By the convexity of $P$, since $m=4>2$, we have ${\bf x}^*\neq {\bf x}^{**}$. Consequently, we get
$$\varphi({\bf v})\geq 2$$
and hence
$$\tau=\varpi({\bf v})+\varphi({\bf v})\geq 3+2=5.\eqno(19)$$

\begin{figure}[h!]
\centering
\includegraphics[scale=.45]{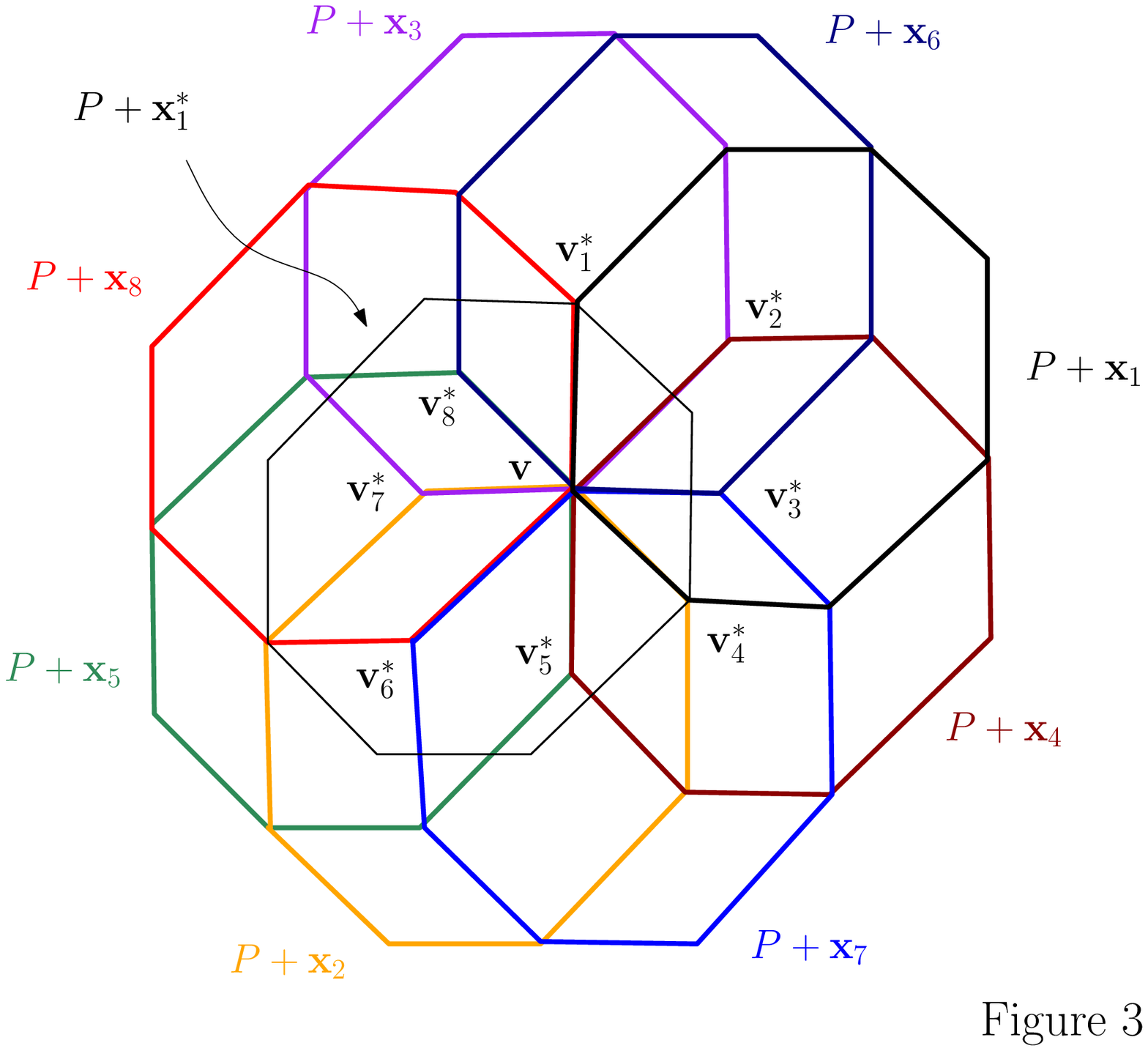}
\end{figure}
If $\ell=0$, then one can deduce that $P+X^{\bf v}$ is an adjacent wheel of eight translates $P+{\bf x}_1,$ $P+{\bf x}_2,$ $\ldots,$ $P+{\bf x}_8$ as shown in Figure 3. Let ${\bf v}_1^*,$ ${\bf v}_2^*,$ ${\bf v}_3^*,$ ${\bf v}_4^*,$ ${\bf v}_5^*,$ ${\bf v}_6^*,$ ${\bf v}_7^*,$ ${\bf v}_8^*$ be the corresponding points of ${\bf v}$ on $F_1+{\bf x}_1,$ $F_2+{\bf x}_4,$ $F_3+{\bf x}_7,$ $F_4+{\bf x}_2,$ $F_5+{\bf x}_5,$ $F_6+{\bf x}_8,$ $F_7+{\bf x}_3,$ $F_8+{\bf x}_6$, respectively. Since ${\bf v}$ is a proper point on $G$, for all $1\leq i\leq 8$ we have ${\bf v}_i^*\not\in K(G)+X$. By Lemma 7, one can deduce that, for each ${\bf v}_i^*$, there is ${\bf x}^*_i\in X$ satisfying both
$${\bf v}_i^*\in B(G)+{\bf x}^*_i$$
and
$${\bf v}\in int (P)+{\bf x}^*_i.$$
If ${\bf x}^*_1={\bf x}^*_2=\cdots={\bf x}^*_8$, we project $P+X^{\bf v}$ and $P+{\bf x}^*_1$ onto the $xy$-plane, as shown by Figure 3. Without confusion, we also use $P+{\bf x}_i$ to denote its projection. By convexity, it is known that $(P+{\bf x}_1)\cap (P+{\bf x}_1^*)$ is centrally symmetric. Furthermore, since they intersect at ${\bf v}^*_1$ and ${\bf v}^*_4$, their intersection must be a parallelogram. Similarly, $(P+{\bf x}_8)\cap (P+{\bf x}_1^*)$ is also a parallelogram. Therefore, the projection of $P+{\bf x}_1^*$ must be an hexagon, which contradicts the assumption that $m=4$. Then we get
$$\varphi({\bf v})\geq 2$$
and
$$\tau=\varpi({\bf v})+\varphi({\bf v})\geq 3+2=5.\eqno(20)$$
		
	
\medskip\noindent
{\bf Case 3.} {\it $\varpi({\bf v})=2$ for all proper points ${\bf v}$}. Then one can deduce that $P+X^{\bf v}$ is an adjacent wheel of five translates $P+{\bf x}_1,$ $P+{\bf x}_2,$ $\ldots,$ $P+{\bf x}_5$. By re-enumeration, we may assume that $P+{\bf x}_1,$ $P+{\bf x}_2,$ $P+{\bf x}_3$ and $P+{\bf x}_4$ are E-type pieces at ${\bf v}$, and $P+{\bf x}_5$ is an F-type piece at ${\bf v}$. For convenience, we may assume that  ${\bf v}\in rint (G_1)+{\bf x}_1$ (as shown in Figure 4).
		
\begin{figure}[h!]
\centering
\includegraphics[scale=.45]{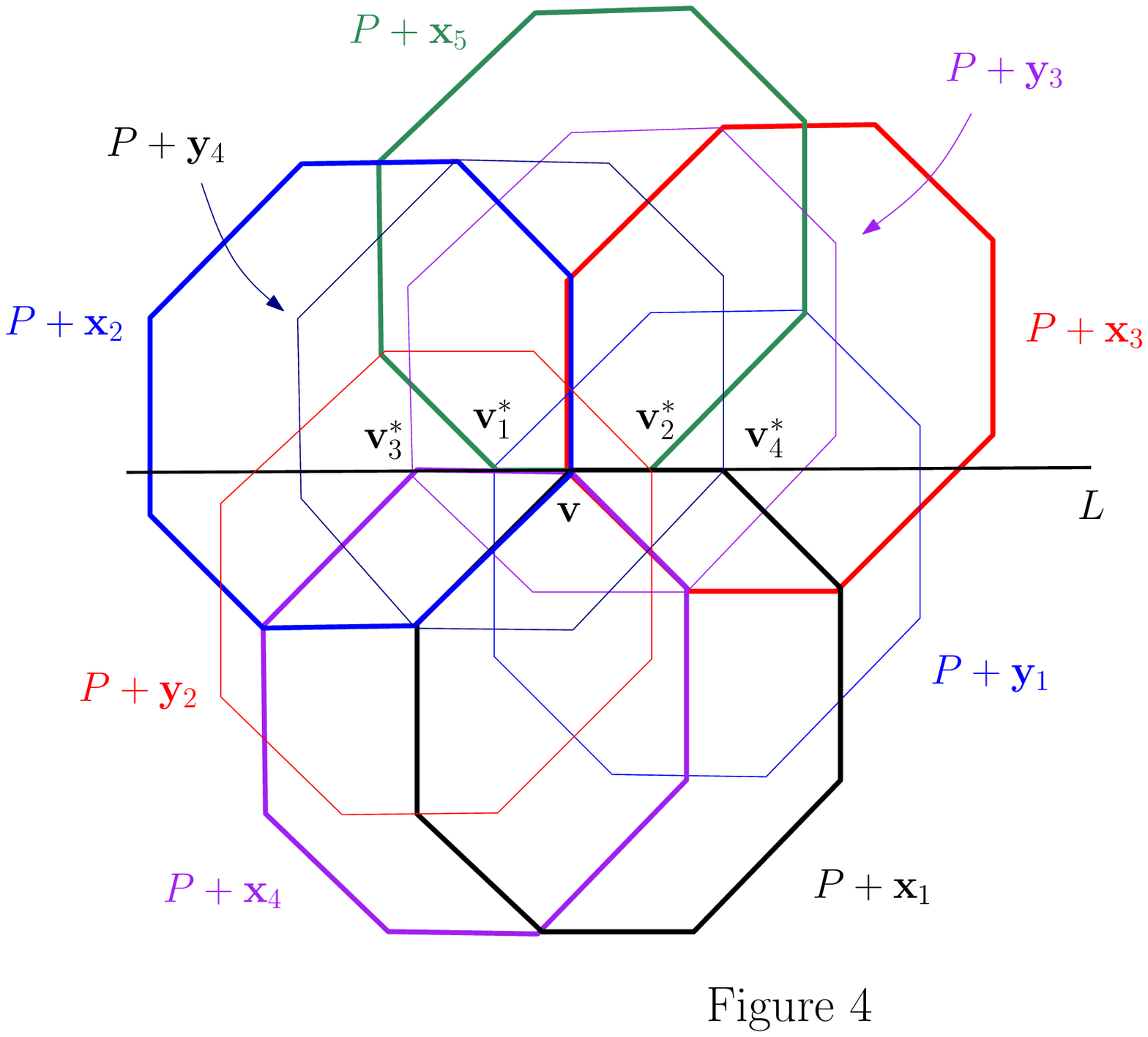}
\end{figure}

As shown in Figure 4, we define ${\bf v}_3^*={\bf v}-{\bf g}_1$ and ${\bf v}_4^*={\bf v}+{\bf g}_1$. Let ${\bf v}_1^*$ be the corresponding point of ${\bf v}$ in $(F_1+{\bf x}_4)\cap (F_5+{\bf x}_5)$ and let ${\bf v}_2^*$ be the corresponding point of ${\bf v}$ in $(F_1+{\bf x}_1)\cap (F_5+{\bf x}_5)$. Then we have ${\bf v}\in {\bf v}_1^*+R_1$ and ${\bf v}\in {\bf v}_2^*-R_1$. By Lemma 9, we may assume that all ${\bf v}_1^*,$ ${\bf v}_2^*,$ ${\bf v}_3^*$ and ${\bf v}_4^*$ are proper points in the respective edges. Therefore, we have
$$\varpi({\bf v}_1^*)=\varpi({\bf v}_2^*)=\varpi({\bf v}_3^*)=\varpi({\bf v}_4^*)=2,\eqno(21)$$
which implies that the adjacent wheels at ${\bf v}_i^*$, where $i=1,2,3,4$, are the same type. Furthermore, since Figure 4 is the only configuration type for $\varpi({\bf v})=2$, all these adjacent wheels have four E-type pieces and one F-type piece. Noting that $P+{\bf x}_4$ is the only F-type piece at ${\bf v}_1^*$, one can deduce that there is ${\bf y}_1\in X$ such that ${\bf v}_1^*\in G_8+{\bf y}_1$. Recall that ${\bf v}\in {\bf v}_1^*+R_1$, by Lemma 5, we have that $({\bf v}_1^*+R_1)\setminus \{{\bf v}_1^*\}\subset  int (P)+{\bf y}_1$ and hence
$${\bf v}\in  int (P)+{\bf y}_1.$$
Furthermore, there is ${\bf z}_1\in X$ such that ${\bf v}_1^*\in G_5+{\bf z}_1$. Since ${\bf v}_3^*= {\bf v}-{\bf g}_1$, from Lemma 6 and the definition of corresponding points, one can deduce that ${\bf v}_3^*={\bf v}_1^*+{\bf g}$ where ${\bf g}\in -R_1$, i.e., ${\bf v}_3^*\in {\bf v}_1^*-R_1\subset P+{\bf z}_1$. Hence, $P+{\bf z}_1$ is in the adjacent wheel at ${\bf v}_3^*$. It is clear that $P+{\bf z}_1$ is an F-type piece at ${\bf v}_3^*$. It follows that there is ${\bf y}_3\in X$ such that ${\bf v}_3^*\in G_7+{\bf y}_3$. From Lemma 5, we know that
$${\bf v}={\bf v}_3^*+{\bf g}_1\in  int (P)+{\bf y}_3.$$
Similarly, there are ${\bf y}_2,$ ${\bf y}_4\in X$ such that ${\bf v}_2^*\in G_3+{\bf y}_2,$ ${\bf v}_4^*\in G_4+{\bf y}_4$,
$${\bf v}\in int (P)+{\bf y}_2$$
and
$${\bf v}\in int (P)+{\bf y}_4.$$
We project $P+X^{\bf v}$ and the corresponding $P+{\bf y}_i$ onto the $xy$-plane and still use $P+{\bf x}_i$ and $P+{\bf y}_i$ to denote their projections. Clearly, we have ${\bf y}_1\not= {\bf y}_2$, ${\bf y}_1\not= {\bf y}_3$ and ${\bf y}_2\not= {\bf y}_4$. Similar to Subcase 3.1.2 of Lemma 3.8 of Yang and Zong \cite{yz2}, one can deduce that ${\bf y}_1= {\bf y}_4$ and ${\bf y}_2 = {\bf y}_3$ can not hold simultaneously. Thus, we have
$$\varphi({\bf v})\geq 3$$
and
$$\tau=\varpi({\bf v})+\varphi({\bf v})\geq 2+3=5.\eqno(22)$$

The lemma is proved. \hfill{$\Box$}

\medskip
\noindent
{\bf Lemma 12.} {\it If $m=5$, then $\tau\geq 5$.}

\medskip\noindent
{\bf Proof.} Assume that $G'\in \mathcal{G}+X$ and let ${\bf v}$ be a proper point in $G'$. By Lemma 8, we have
$$\varpi({\bf v})=\kappa\cdot2+\ell\cdot\frac12\geq 2.\eqno(23)$$
We consider the following three cases.
	
\medskip\noindent
{\bf Case 1.}  {\it $\varpi({\bf v})\geq 4$ for some proper point ${\bf v}$ in $G'$.} By Lemma 10, one obtains
$$\varphi({\bf v})\geq\left\lceil{\frac{m-3}{2}}\right\rceil = 1,$$
and hence
$$\tau=\varpi({\bf v})+\varphi({\bf v})\geq 4+1=5.\eqno(24)$$
	
\medskip\noindent
{\bf Case 2.} {\it $\varpi({\bf v})=3$ for some proper point ${\bf v}$ in $G'$.} Since
$$\varpi({\bf v})=\kappa\cdot2+\ell\cdot\frac12=3$$
and $\kappa$ is a positive integer, we get $\ell>0$. Applying the arguments used in Lemma 11 (the first part of Case 2), one can obtain $$\varphi({\bf v})\geq 2,$$
and hence
$$\tau=\varpi({\bf v})+\varphi({\bf v})\geq 3+2=5.\eqno(25)$$
	
\medskip\noindent
{\bf Case 3.} {\it $\varpi({\bf v})=2$ for all proper points ${\bf v}\in \mathcal{G}+X$.} Then we have $\kappa=1$ and $\ell=0$. One can deduce that $P+X^{\bf v}$ is an adjacent wheel of five translates $P+{\bf x}_1,$ $P+{\bf x}_2,$ $\ldots,$ $P+{\bf x}_5$, all of which are E-type at ${\bf v}$. By re-enumeration, we may assume that $P+{\bf x}_1,$ $P+{\bf x}_2,$ $\ldots,$ $P+{\bf x}_5$ are arranged in the order as shown in Figure 5. Let ${\bf v}_1^*,$ ${\bf v}_2^*,$ ${\bf v}_3^*,$ ${\bf v}_4^*,$ ${\bf v}_5^*$ be the corresponding points of ${\bf v}$ in $F_1+{\bf x}_1,$ $F_3+{\bf x}_4,$ $F_5+{\bf x}_2,$ $F_7+{\bf x}_5,$ $F_9+{\bf x}_3$, respectively. By Lemma 7, for each of these five points ${\bf v}_i^*$, there is ${\bf y}_i\in X$ satisfying both
$${\bf v}_i^*\in B(G)+{\bf y}_i$$
and
$${\bf v}\in int (P)+{\bf y}_i.$$

\begin{figure}[h!]
\centering
\includegraphics[scale=.45]{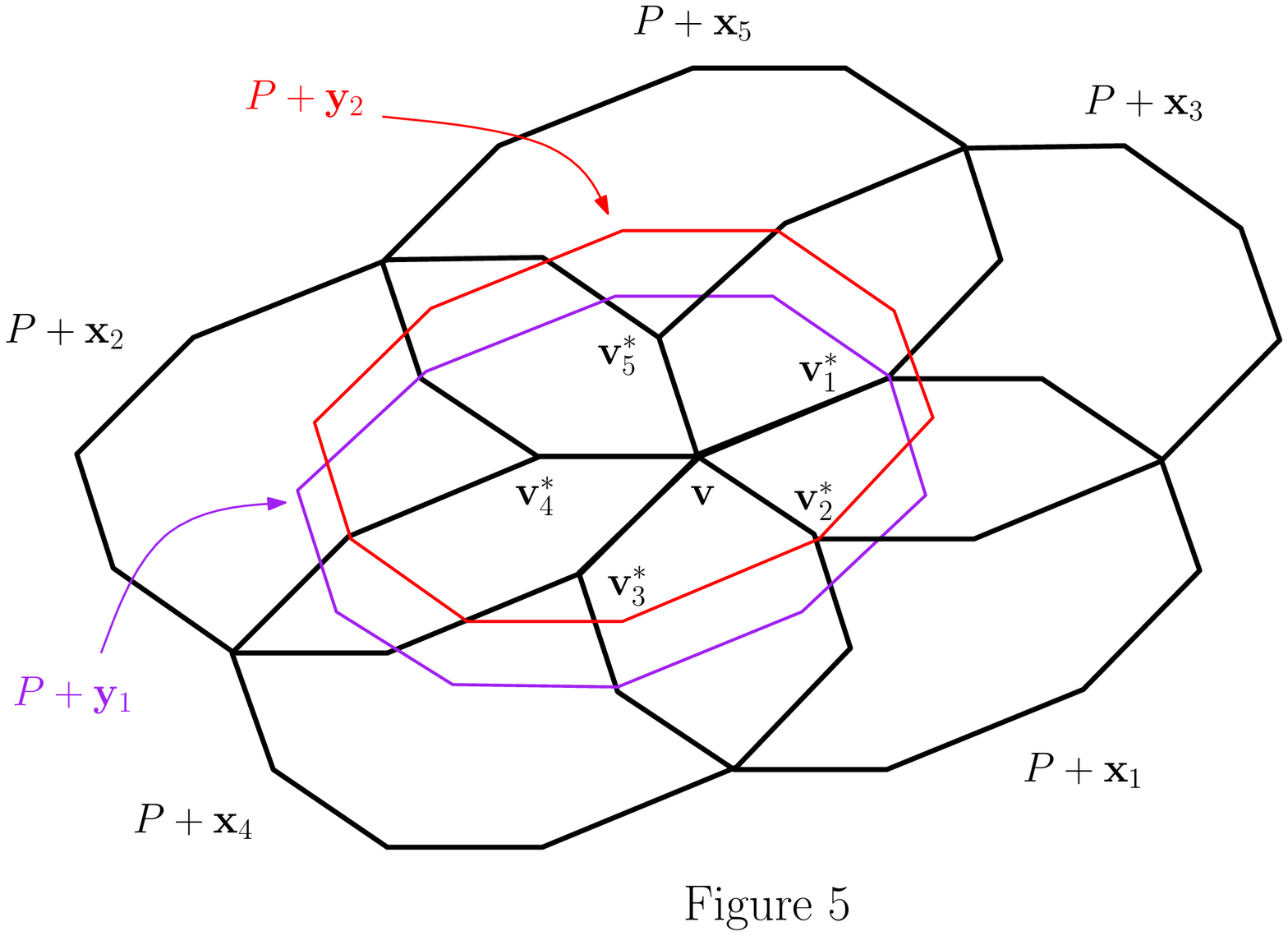}
\end{figure}

We project $P+X^{\bf v}$ and $P+{\bf y}_i$ onto the $xy$-plane and denote their projections still by $P+{\bf x}_i$ and $P+{\bf y}_i$, respectively. Since $\varpi({\bf v})=2$ holds for all proper points ${\bf v}\in \mathcal{G}+X$, the adjacent wheel at each of the five points ${\bf v}_i^*$ is unique. If ${\bf y}_1={\bf y}_2={\bf y}_3,$ one can deduce that $(P+{\bf y}_1)\cap (P+{\bf x}_1)$ is both a parallelogram and an octagon, which is impossible. If ${\bf y}_1={\bf y}_2={\bf y}_4,$ both $(P+{\bf y}_1)\cap (P+{\bf x}_4)$ and $(P+{\bf y}_1)\cap (P+{\bf x}_5)$ are parallelograms and therefore the projection of $P$ is an hexagon, which contradicts the assumption that $m=5$. Consequently, one can deduce that no three of the five $P+{\bf y}_i$ are identical. Therefore, we have
$$\varphi({\bf v})\geq 3$$
and hence
$$\tau=\varpi({\bf v})+\varphi({\bf v})\ge 2+3=5.\eqno(26)$$

The lemma is proved. \hfill{$\Box$}

\medskip
\noindent
{\bf Lemma 13.} {\it If $m=6$, then $\tau\geq 6$.}

\medskip\noindent
{\bf Proof.} Assume that $G'\in \mathcal{G}+X$ and let ${\bf v}$ be a proper point in $G'$. By Lemma 8, we get
$$\varpi({\bf v})=\kappa\cdot\frac52+\ell\cdot\frac12\geq 3.\eqno(27)$$
Now, we consider the following two cases.
		
\medskip
\noindent
{\bf Case 1.} {\it $\varpi({\bf v})\geq 4$ for some proper point ${\bf v}$ in $G'$.} By Lemma 10, we have
$$\varphi({\bf v})\geq\left\lceil{\frac{m-3}{2}}\right\rceil = 2$$
and hence
$$\tau=\varpi({\bf v})+\varphi({\bf v})\geq 4+2=6.\eqno(28)$$
		
\medskip
\noindent
{\bf Case 2.} {\it $\varpi({\bf v})=3$ for all proper point ${\bf v}$ in $G'$.} Then, we have $\ell=1$. In other words, there is an F-type piece $P+{\bf x}$ at ${\bf v}$, where ${\bf x}\in X$. Suppose that ${\bf v}\in rint (F_i)+{\bf x}$. By  Lemma 9, one can choose ${\bf v}^*\in (G_i+{\bf x})\setminus(K(G)+X)$ and ${\bf v}^{**}={\bf v}^*+{\bf g}_i\in (G_{i+1}+x)\setminus(K(G)+X)$ such that
$${\bf v}\in {\bf v}^*+R_i={\bf v}^{**}-R_i.$$
By applying Lemma 7 to the point ${\bf v}^*$, there are at least two different points ${\bf y}_1^*$ and ${\bf y}_2^*$ in $X$ such that
$${\bf v}^*\in B(G)+{\bf y}_1^*,$$
$${\bf v}^*\in B(G)+{\bf y}_2^*,$$
$$({\bf v}^*+R_{i})\setminus\{{\bf v}^*\}\subset  int(P)+{\bf y}_1^*$$
and
$$({\bf v}^*+R_i)\setminus\{{\bf v}^*\}\subset  int(P)+{\bf y}_2^{*}.$$
Similarly, by applying Lemma 7 to the point ${\bf v}^{**}$, there are at least two different points ${\bf y}_1^{**}$ and ${\bf y}_2^{**}$ such that
$${\bf v}^{**}\in B(G)+{\bf y}_1^{**},$$
$${\bf v}^{**}\in B(G)+{\bf y}_2^{**},$$
$$({\bf v}^{**}-R_{i})\setminus\{{\bf v}^{**}\}\subset int(P)+{\bf y}_1^{**}$$
and
$$({\bf v}^{**}-R_{i})\setminus\{{\bf v}^{**}\}\subset  int(P)+{\bf y}_2^{**}.$$
By convexity, it is easy to see that ${\bf y}_1^*,$ ${\bf y}_2^*,$ ${\bf y}_1^{**}$ and ${\bf y}_2^{**}$ are pairwise distinct. Therefore, we get
$$\varphi({\bf v})\geq 4$$
and
$$\tau=\varpi({\bf v})+\varphi({\bf v})\geq 3+4=7.\eqno(29)$$

The lemma is proved. \hfill{$\Box$}

\medskip
\noindent
{\bf Lemma 14.} {\it If $m=7$, then $\tau\geq 6$.}

\medskip\noindent
{\bf Proof.} Assume that $G'\in \mathcal{G}+X$ and let ${\bf v}$ be a proper point on $G'$. From Lemma 8, we have
$$\varpi({\bf v})=\kappa\cdot3+\ell\cdot\frac12\geq 3.\eqno(30)$$
Now, we consider the following two cases.
	
\medskip\noindent
{\bf Case 1.} {\it $\varpi({\bf v})\geq 4$ for some proper point ${\bf v}$ in $G'$.} By Lemma 10, we have
$$\varphi({\bf v})\geq\left\lceil{\frac{m-3}{2}}\right\rceil = 2$$
and hence
$$\tau=\varpi({\bf v})+\varphi({\bf v})\geq 4+2=6.\eqno(31)$$
	
\medskip\noindent
{\bf Case 2.} {\it $\varpi({\bf v})=3$ for all proper points ${\bf v}$ in $\mathcal{G}+X$.} Then, by (30) we have $\ell =0$. Consequently, $P+X^{\bf v}$ is an adjacent wheel of seven translates $P+{\bf x}_1,$ $P+{\bf x}_2,$ $\ldots,$ $P+{\bf x}_7$, all of them are E-type pieces at ${\bf v}$.  By re-enumeration, we may assume that $P+{\bf x}_1,$ $P+{\bf x}_2,$ $\ldots,$ $P+{\bf x}_7$ are arranged in the order as shown in Figure 6.  Let ${\bf v}_1^*,$ ${\bf v}_2^*,$ ${\bf v}_3^*,$ ${\bf v}_4^*,$ ${\bf v}_5^*,$ ${\bf v}_6^*$ and ${\bf v}_7^*$ be the corresponding points of ${\bf v}$ in $F_1+{\bf x}_1,$ $F_3+{\bf x}_6,$ $F_5+{\bf x}_4,$ $F_7+{\bf x}_2,$ $F_9+{\bf x}_7,$ $F_{11}+{\bf x}_5$ and $F_{13}+{\bf x}_3$, respectively. By Lemma 7, for each of these seven points ${\bf v}_i^*$, there are two different points ${\bf y}_i, {\bf y}_i^*\in X$ satisfying
$${\bf v}_i^*\in B(G)+{\bf y}_i,$$
$${\bf v}_i^*\in B(G)+{\bf y}^*_i,$$
$${\bf v}\in int (P)+{\bf y}_i$$
and
$${\bf v}\in int (P)+{\bf y}^*_i.$$

\begin{figure}[h!]
\centering
\includegraphics[scale=.45]{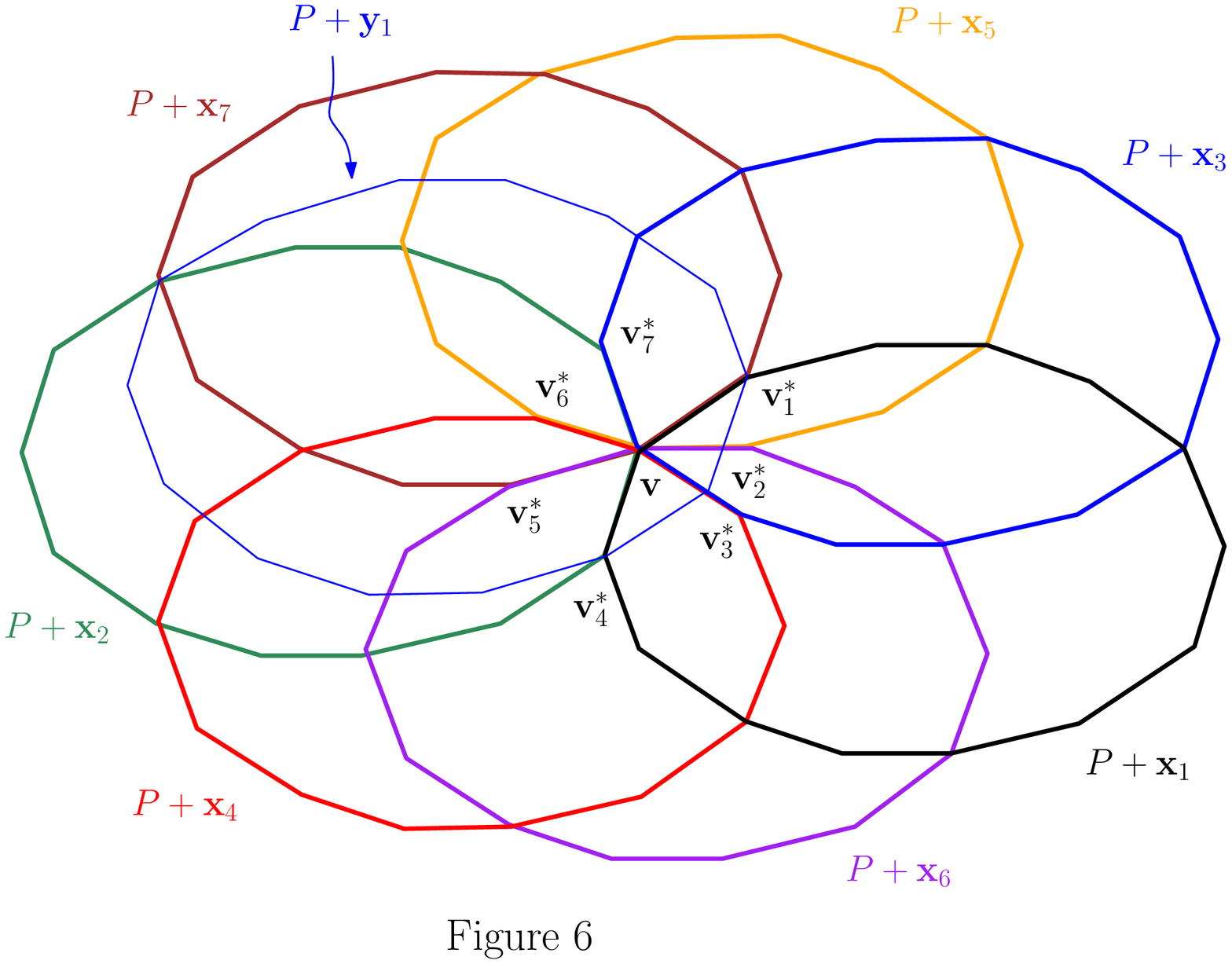}
\end{figure}

We project $P+X^{\bf v}$ and the corresponding $P+{\bf y}_i$ onto the $xy$-plane and denote the projections still by $P+{\bf x}_i$ and $P+{\bf y}_i$, respectively. If all ${\bf v}_i^*$ are vertices of $P+{\bf y}_1$, by convexity and symmetry, it is easy to see that both $(P+{\bf x}_1)\cap (P+{\bf y}_1)$ and $(P+{\bf x}_7)\cap (P+{\bf y}_1)$ are parallelograms. Consequently, we get that $P+{\bf y}_1$ is an hexagon, which contradicts the assumption $m=7$. If, without loss of generality, ${\bf v}^*_3$ is not a vertex of $P+{\bf y}_1$, then ${\bf y}_1$, ${\bf y}_3$ and ${\bf y}_3^*$ are pairwise distinct. Consequently, we get
$$\varphi({\bf v})\geq 3$$
and therefore
$$\tau=\varpi({\bf v})+\varphi({\bf v})\geq 3+3=6.\eqno(32)$$

The lemma is proved. \hfill{$\Box$}

\medskip
\noindent
{\bf Lemma 15.} {\it If $m\geq 8$, then $\tau\geq 7$.}

\medskip\noindent
{\bf Proof.} Assume that $G'\in \mathcal{G}+X$ and let ${\bf v}$ be a proper point on $G'$. By Lemma 8 and Lemma 10, we have
$$\varpi({\bf v})=\kappa\cdot\frac{m-1}{2}+\ell\cdot\frac12\geq 4,\eqno(33)$$
$$\varphi({\bf v})\geq \left\lceil{\frac{m-3}{2}}\right\rceil\geq 3\eqno(34)$$
and
$$\tau=\varpi({\bf v})+\varphi({\bf v})\geq 4+3=7.\eqno(35)$$
The lemma is proved. \hfill{$\Box$}

\bigskip
\noindent
{\bf Proof of Theorem 1.} If $P+X$ is a threefold tiling of $\mathbb{E}^3$, it follows from Lemmas 11, 12, 13, 14 and 15 that every belt of $P$ has at most six facets. Then, it follows by Lemma 2 that $P$ must be a parallelohedron. Furthermore, by Lemma 1, $P$ must be a parallelotope, a hexagonal prism, a rhombic dodecahedron, an elongated dodecahedron, or a truncated octahedron. The theorem is proved. \hfill{$\Box$}

\medskip
\noindent
{\bf Proof of Theorem 2.} If $P+X$ is a fourfold tiling of $\mathbb{E}^3$, it follows from Lemmas 11, 12, 13, 14 and 15 that every belt of $P$ has at most six facets. Then, it follows by Lemma 2 that $P$ must be a parallelohedron. Furthermore, by Lemma 1, $P$ must be a parallelotope, a hexagonal prism, a rhombic dodecahedron, an elongated dodecahedron, or a truncated octahedron. The theorem is proved. \hfill{$\Box$}

\medskip
\noindent
{\bf Proof of Theorem 3.} If $P+X$ is a fivefold tiling of $\mathbb{E}^3$, it follows from Lemmas 11, 12, 13, 14 and 15 that every belt of $P$ has at most ten facets. The theorem is proved. \hfill{$\Box$}

\medskip
\noindent
{\bf Proof of Theorem 4.} If $P+X$ is a sixfold tiling of $\mathbb{E}^3$, it follows from Lemmas 11, 12, 13, 14 and 15 that every belt of $P$ has at most fourteen facets. The theorem is proved. \hfill{$\Box$}

\bigskip
To end this paper, we emphasize the following open problems.

\medskip\noindent
{\bf Problem 1.} To characterize all the three-dimensional fivefold lattice tiles.

\medskip\noindent
{\bf Problem 2.} To characterize all the three-dimensional fivefold translative tiles.

\vspace{0.6cm}\noindent
{\bf Acknowledgements.} This work is supported by ERC Starting Grant No. 713927, the National Natural Science Foundation of China (NSFC11921001), the National Key Research and Development Program of China (2018YFA0704701), and 973 Program 2013CB834201.

\vspace{0.6cm}
\noindent
Mei Han, Center for Applied Mathematics, Tianjin University, Tianjin 300072, China.

\noindent
Qi Yang, Department of Mathematics, Bar Ilan University, Israel.

\noindent
Kirati Sriamorn, Department of Mathematics and Computer Science, Chulalongkorn University, Thailand.

\noindent
Chuanming Zong, Center for Applied Mathematics, Tianjin University, Tianjin 300072, China.

\noindent
Email: cmzong@math.pku.edu.cn

\begin{thebibliography}{99}
	\bibitem{alek}A. D. Aleksandrov, On tiling space by polytopes, {\it Vestnik Leningrad Univ. Ser. Mat. Fiz. Him}. {\bf 9} (1954), 33-43.
	\bibitem{boll}U. Bolle, On multiple tiles in $R^2$, {\it Intuitive Geometry,} Colloq. Math. Soc. J. Bolyai {\bf 63} (1994), 39-43.
	\bibitem{delo}B. N. Delone, Sur la partition reguli$\grave{e}$re de l'espace $\grave{a}$ $4$ dimensions I, II, {\it Izv. Akad. Nauk SSSR, Ser. VII} (1929), 79-110; 147-164.
	\bibitem{dgsw}M. Dutour Sikiri$\acute{\rm c}$, A. Garber, A. Sch$\ddot{\rm u}$rmann and C. Waldmann, The complete classification of five-dimensional Dirichlet-Voronoi polyhedra of translational lattices, {\it Acta Crystallogr. Sect. A}  {\bf 72} (2016), 673-683.
	\bibitem{enge}P. Engel, On the symmetry classification of the four-dimensional parallelohedra, {\it Z. Kristallographie} {\bf 200} (1992), 199-213.
	\bibitem{fedo}E. S. Fedorov, Elements of the study of figures, {\it Zap. Mineral. Imper. S. Petersburgskogo Ob$\check{s}$$\check{c}$}, {\bf 21}(2) (1885), 1-279.
	\bibitem{furt}P. Furtw\"angler, \"Uber Gitter konstanter Dichte, {\it Monatsh. Math. Phys.} {\bf 43} (1936), 281-288.
	\bibitem{grs}N. Gravin, S. Robins and D. Shiryaev, Translational tilings by a polytope, with multiplicity. {\it Combinatorica} {\bf 32} (2012), 629-649.
	\bibitem{gkrs}N. Gravin, M. N. Kolountzakis, S. Robins and D. Shiryaev, Structure results for multiple tilings in 3D. {\it Discrete Comput. Geom.} {\bf 50} (2013), 1033-1050.
	\bibitem{Grepstad}S. Grepstad and N. Lev, Multi-tiling and Riesz bases. {\it Adv. Math.}  {\bf 252} (2014), 1-6.
	\bibitem{hajo}G. Haj\'os, \"Uber einfache und mehrfache Bedeckung des $n$-dimensionalen Raumes mit einem W\"urfelgitter, {\it Math. Z.}
	{\bf 47} (1941), 427-467.
    \bibitem{hsyz}M. Han, K. Sriamorn, Q. Yang and C. Zong, The twofold translative tiles in three-dimensional space, arXiv:2106.15388.
	\bibitem{kolo}M. N. Kolountzakis, Multiple lattice tiles and Riesz bases of exponentials. {\it Proc. Amer. Math. Soc.} {\bf 143} (2015), 741-747.
	\bibitem{Lev-Liu}N. Lev and B. Liu, Multi-tiling and equidecomposability of polytopes by lattice translates. {\it Bull. Lond. Math. Soc.}  {\bf 51}  (2019), 1079-1098.
	\bibitem{mcmu}P. McMullen, Convex bodies which tiles space by translation, {\it Mathematika} {\bf 27} (1980), 113-121.
	\bibitem{mink}H. Minkowski, Allgemeine Lehrs\"atze \"uber konvexen Polyeder, {\it Nachr. K. Ges. Wiss. G\"ottingen, Math.-Phys. KL}. (1897), 198-219.
	\bibitem{robi}R. M. Robinson, Multiple tilings of $n$-dimensional space by unit cubes, {\it Math. Z.} {\bf 166}
	(1979), 225-275.
	\bibitem{sriamorn}K. Sriamorn, Twofold translative tilings with convex bodies, arXiv:1601.04312.
	\bibitem{stog}M. I. $\check{S}$togrin, Regular Dirichlet-Voronoi partitions for the second triclinic group (in Russian), {\it Proc. Steklov. Inst. Math.} {\bf 123} (1975).
	\bibitem{venk}B. A. Venkov, On a class of Euclidean polytopes, {\it Vestnik Leningrad Univ. Ser. Mat. Fiz. Him}. {\bf 9} (1954), 11-31.
	\bibitem{voro}G. F. Voronoi,  Nouvelles applications des paramm\`etres continus \`a la th\'eorie des formes quadratiques.
	Deuxi\`eme M\'emoire. Recherches sur les parall\'elo\`edres primitifs, {\it J. reine angew. Math.} {\bf 134} (1908), 198-287; {\bf 135} (1909), 67-181.
	\bibitem{yz1}Q. Yang and C. Zong, Multiple lattice tiling in Euclidean spaces, {\it Canadian Math. Bull.} {\bf 62} (2019), 923-929.
	\bibitem{yz2}Q. Yang and C. Zong, Characterization of the two-dimensional fivefold translative tiles, {\it Bull. Soc. Math. France} {\bf 149} (2021), 119-153.
	\bibitem{zong96}C. Zong, {\it Strange Phenomena in Convex and Discrete Geometry.} Springer-Verlag, New York, 1996.
	\bibitem{zong05}C. Zong, What is known about unit cubes. {\it Bull. Amer. Math. Soc.} {\bf 42} (2005), 181-211.
	\bibitem{zong06}C. Zong, {\it The Cube: A Window to Convex and Discrete Geometry.} Cambridge University Press, Cambridge, 2006.
	\bibitem{zong20}C. Zong, Can you pave the plane with identical tiles? {\it Notices Amer. Math. Soc.} {\bf 67} (2020), 635-646.
	\bibitem{zong-x}C. Zong, Characterization of the two-dimensional five- and sixfold lattice tiles, arXiv:1712.01122, arXiv:1904.06911.
\end{thebibliography}
\end{document}